%% file: ICBS_Wedrich-arxiv.tex
\documentclass[10pt,twoside]{article}

\newcommand{\email}[1]{{\footnotesize {\it E-mail address}: {\tt #1\vphantom{y}}}}

\usepackage[english]{babel}
\usepackage[utf8]{inputenc}
\usepackage{multicol}
\usepackage{tabu}
\usepackage{amsmath}
\usepackage{amssymb}
\usepackage{enumitem}
\setlist{itemsep=1pt, topsep=4pt, parsep=0pt, partopsep=0pt}
\usepackage{bbm}
\usepackage{colortbl}
\usepackage{hyperref}

\makeatletter
\let\oldthebibliography\thebibliography
\let\endoldthebibliography\endthebibliography
\renewenvironment{thebibliography}[1]{%
  \oldthebibliography{#1}%
  \setlength{\itemsep}{-4pt}%
  \setlength{\parsep}{0pt}%
  \setlength{\topsep}{-5pt}%
  \setlength{\partopsep}{0pt}%
}{\endoldthebibliography}
\makeatother

\numberwithin{equation}{section}

\newtheorem{thm}{Theorem}[section]
\newtheorem{theorem}[thm]{Theorem}
\newtheorem{exa}[thm]{Example}

\newtheorem{definition}[thm]{Definition}

\newtheorem{conjecture}[thm]{Conjecture}

\newtheorem{rem}[thm]{Remark}

\usepackage{graphicx}
\usepackage{xcolor}
\newcommand{\arxiv}[1]{\href{https://arxiv.org/abs/#1}{\small  arXiv:#1}}

\newcommand{\googlebooks}[1]{(preview at \href{https://books.google.com/books?id=#1}{google books})}

\newcommand{\numdam}[1]{}

\def\Sz{\mathcal{S}_0^N}
\def\Sztwo{\mathcal{S}_0^2}
\def\KhRN{\operatorname{KhR}_N}
\def\CKhRN{\operatorname{CKhR}_N}
\def\gl{\mathfrak{gl}}
\def\glN{\gl_N}
\def\GLN{\mathrm{GL}(N)}
\def\Z {{\mathbb{Z}}}

\def\R {{\mathbb{R}}} 

\def\Q {{\mathbb{Q}}}
\def\kk {{\mathbbm{k}}}

\newcommand{\fourm}{4}

\newcommand{\target}{\mathcal{V}}

\newcommand{\Links}{\mathbf{Links}}
\newcommand{\Kb}{\mathrm{K}^b}
\newcommand{\Der}{\mathrm{D}}

\newcommand{\gmod}{{-}\mathrm{gmod}}
\newcommand{\gpmod}{{-}\mathrm{gpmod}}
\newcommand{\linkhom}{\mathrm{H}}
\newcommand{\bigrabgrps}{\mathrm{gr}^{\Z\times \Z}\mathbf{AbGrp}}

\input{figures.tex}

\begin{document}

 \markboth{\hfill{\rm Paul Wedrich} \hfill}{\hfill {\rm From Link Homology to Topological Quantum Field Theories \hfill}}

\title{From Link Homology to Topological Quantum Field Theories}

\author{Paul Wedrich}

\date{}

\maketitle

\begin{abstract}
This survey reviews recent advances connecting link homology theories to
invariants of smooth 4-manifolds and extended topological quantum field
theories. Starting from joint work with Morrison and Walker, I explain how
functorial link homologies that satisfy additional invariance conditions become
diagram-independent, give rise to braided monoidal 2-categories, extend
naturally to links in the 3-sphere, and globalize to skein modules for
4-manifolds. Later developments show that these skein lasagna modules furnish
invariants of embedded and immersed surfaces and admit computation via handle
decompositions. I then survey structural properties, explicit computations, and
applications to exotic phenomena in 4-manifold topology, and place link homology
and skein lasagna modules within the framework of extended topological quantum
field theories.
\end{abstract}

\input{maintext}





\def\address#1{\begingroup\setlength{\parindent}{0pt}%
\vspace*{15pt}
\fontsize{10}{8}\selectfont\sc #1
\endgroup}

 \address{Fachbereich Mathematik, Universit\"at Hamburg, 
 Bundesstra{\ss}e 55, 
 20146 Hamburg, Germany
 \href{https://paul.wedrich.at}{paul.wedrich.at}}
 \email{paul.wedrich@uni-hamburg.de}

\end{document}

%% file: figures.tex
\usepackage{tikz}
\usepackage{tikz-3dplot}
\usetikzlibrary{3d, decorations.pathmorphing}
\tikzset{anchorbase/.style={baseline={([yshift=-0.5ex]current bounding box.center)}}}
\usetikzlibrary{calc}
\usetikzlibrary{decorations.markings}
\usetikzlibrary{decorations.pathreplacing}
\usetikzlibrary{arrows,shapes,positioning}
\tikzstyle directed=[postaction={decorate,decoration={markings,
    mark=at position #1 with {\arrow{>}}}}]
\tikzstyle rdirected=[postaction={decorate,decoration={markings,
    mark=at position #1 with {\arrow{<}}}}]
    \tikzset{
    tinynodes/.style={font=\tiny,text height=0.75ex,text depth=0.15ex},
    smallnodes/.style={font=\scriptsize,text height=0.75ex,text depth=0.15ex}
  }

\newcommand{\fsphere}[3]{
\draw[white, line width=1mm] (#1-#3,#2) to [out=270,in=180] (#1,#2-#3*2/3) to [out=0,in=270] (#1+#3,#2);
\draw (#1-#3,#2) to [out=270,in=180] (#1,#2-#3*2/3) to [out=0,in=270] (#1+#3,#2);
\draw[thick] (#1,#2) circle (#3);
}

\newcommand{\bsphere}[3]{
\draw[dashed] (#1-#3,#2) to [out=90,in=180] (#1,#2+#3*2/3) to [out=0,in=90] (#1+#3,#2);
}

\newcommand{\sphere}[3]{
\bsphere{#1}{#2}{#3}
\fsphere{#1}{#2}{#3}
}

\newcommand{\fgsphere}[3]{
\draw[white, line width=1mm] (#1-#3,#2) to [out=270,in=180] (#1,#2-#3*2/3) to [out=0,in=270] (#1+#3,#2);
\draw[opacity=.3] (#1-#3,#2) to [out=270,in=180] (#1,#2-#3*2/3) to [out=0,in=270] (#1+#3,#2);
\draw[thick, opacity=.3] (#1,#2) circle (#3);
}

\newcommand{\bgsphere}[3]{
\draw[dashed,opacity=.3] (#1-#3,#2) to [out=90,in=180] (#1,#2+#3*2/3) to [out=0,in=90] (#1+#3,#2);
}

\newcommand{\trefoil}[4]{
\begin{scope}[shift={(#1,#2-5.9)},scale=#3]
\draw[#4] (-0.3,5.15) to [out=150,in=270] (-1,6) to   [out=90,in=180] (0,6.8) to [out=0,in=110] (0.85,6.3)
(-0.7,6) to [out=20,in=160] (1,6) to   [out=340,in=60] (1.4,5) to [out=240,in=330] (0.3,4.85)
(0.85,5.7) to [out=240,in=30] (0,5) to   [out=210,in=300] (-1.4,5) to [out=120,in=200] (-1.15,5.85);
\end{scope}
}

\newcommand{\hopflink}[4]{
\begin{scope}[shift={(#1,#2)},scale=#3]
\begin{scope}[yscale=.66]
\draw[#4]
(.5,0) to [out=90,in=0] (-.5,1) to [out=180,in=90] (-1.5,0);
\draw[white, line width=1mm]
(-.5,0) to [out=90,in=180] (.5,1) to [out=0,in=90] (1.5,0);
\draw[#4]
(-.5,0) to [out=90,in=180] (.5,1) to [out=0,in=90] (1.5,0) to [out=270,in=0] (.5,-1) to [out=180,in=270] (-.5,0);
\draw[white, line width=1mm]
(.5,0) to [out=270,in=0] (-.5,-1);
\draw[#4]
(.5,0) to [out=270,in=0] (-.5,-1) to [out=180,in=270] (-1.5,0);
\end{scope}
\end{scope}
}

\newcommand{\unknot}[4]{
\begin{scope}[shift={(#1-1,#2)},scale=#3]
\begin{scope}[yscale=.66]
\draw[#4]
(-.5,0) to [out=90,in=180] (.5,1) to [out=0,in=90] (1.5,0) to [out=270,in=0] (.5,-1) to [out=180,in=270] (-.5,0);
\end{scope}
\end{scope}
}

\newcommand{\torushole}[4]{
\draw[#4] (#1-#3/2,#2+#3/6) to [out=300,in=180] (#1,#2-#3/6) to [out=0,in=240] (#1+#3/2,#2+#3/6);
\draw[#4] (#1-#3/3,#2) to [out=45,in=180] (#1,#2+#3/6) to [out=0,in=135] (#1+#3/3,#2);
}

\newcommand{\torus}[4]{
\draw[#4] (#1,#2) ellipse (#3 and #3*2/3);
\torushole{#1}{#2}{#3}{#4}
}

\newcommand{\lasagnafigure}[1]{
\begin{tikzpicture}[anchorbase,scale=#1]
\bsphere{0}{0}{15}
\sphere{1}{-5}{3}
\sphere{10.5}{-.5}{3}
\sphere{-8}{.6}{3}
\trefoil{-9}{1.5}{1}{red,very thick}
\hopflink{10.5}{0}{1}{red,very thick}
\unknot{0.5}{-5}{.75}{red,very thick}
\unknot{3}{-4.5}{.75}{red,very thick}
\draw[orange, thick]
(-10.5,.7) to [out=80,in=290] (-7.3,11)
(-2.9,11) to [out=280, in=180] (-1,9) to [out=0,in=260] (.7,11)
(5.2,10.7) to [out=250, in=180] (6,8) to [out=0,in=260] (7.3,10)
(10.3,10) to [out=250, in=90] (9,6) to [out=270,in=110] (12,0)
(9,0) to [out=115, in=0] (0,5) to [out=180,in=70] (-7.5,.7)
(3.1,-4.7) to [out=75, in=0] (-.5,3) to [out=180,in=110] (-.8,-5.2)
(1.5,-4.4) to [out=75, in=0] (-.1,.5) to [out=180,in=110] (.7,-5.2)
;
\torushole{-4}{6}{2}{thick, orange}
\torushole{2}{7}{2}{thick, orange}
\torushole{-.3}{1.5}{2}{thick, orange}
\torus{-5}{-4}{2.3}{thick, orange}
\trefoil{3}{9}{1.5}{red,very thick}
\hopflink{-5}{11}{1.5}{red,very thick}
\unknot{9}{10}{1.5}{red,very thick}
\node at (-12,-2.5) {\tiny${\color{red}L_j}$};
\node at (-2,-8.5) {\tiny${\color{red}L_k}$};
\node at (7.5,-3.5) {\tiny${\color{red}L_l}$};
\node at (-2,13) {\tiny${\color{red}L}$};
\node at (-1.5,6) {\tiny${\color{orange}\Sigma}$};
\node at (1.5,1) {\tiny${\color{orange}\Sigma}$};
\fsphere{0}{0}{15}
\end{tikzpicture}
}

\newcommand{\lasagnafillingfigure}[1]{
\begin{tikzpicture}[anchorbase,scale=#1]
    \begin{scope}[shift={(-18,-1)},xscale=3,yscale=12]
        \draw[dashed] (0,1) to [out=0,in=90] (1,0) to [out=270,in=0] (0,-1)
        (0,1) to [out=180,in=90] (-1,0) to [out=270,in=180] (0,-1);
    \end{scope}
    \begin{scope}[shift={(-18,-1)},xscale=9,yscale=12]
        \draw[thick] (0,-1) to [out=15,in=180] (2,-1.3) to [out=0,in=270] (3.8,0)
        to [out=90,in=0] (2,1.3) to [out=180,in=340] (0,1);
    \end{scope}
    \sphere{1}{-5}{3}
    \sphere{10.5}{-.5}{3}
    \sphere{-8}{.6}{3}
    \draw[white, line width=1mm]
    (-9.5,.7) to [out=80,in=290] (-7.3,11)
    (-2.9,11) to [out=280, in=180] (-1,9) to [out=0,in=260] (.7,11)
    (5.2,10.7) to [out=250, in=180] (6,8) to [out=0,in=260] (7.3,10)
    (10.3,10) to [out=250, in=90] (9,6) to [out=270,in=110] (12,0)
    (9,0) to [out=115, in=0] (0,5) to [out=180,in=70] (-6.5,.7)
    (3.1,-4.7) to [out=75, in=0] (-.5,3) to [out=180,in=110] (-.8,-5.2)
    (1.5,-4.4) to [out=75, in=0] (-.1,.5) to [out=180,in=110] (.7,-5.2)
    ;
    \draw[orange, thick]
    (-9.5,.7) to [out=80,in=290] (-7.3,11)
    (-2.9,11) to [out=280, in=180] (-1,9) to [out=0,in=260] (.7,11)
    (5.2,10.7) to [out=250, in=180] (6,8) to [out=0,in=260] (7.3,10)
    (10.3,10) to [out=250, in=90] (9,6) to [out=270,in=110] (12,0)
    (9,0) to [out=115, in=0] (0,5) to [out=180,in=70] (-6.5,.7)
    (3.1,-4.7) to [out=75, in=0] (-.5,3) to [out=180,in=110] (-.8,-5.2)
    (1.5,-4.4) to [out=75, in=0] (-.1,.5) to [out=180,in=110] (.7,-5.2)
    ;
    \torushole{2}{7}{2}{thick, orange}
    \torushole{-.3}{1.5}{2}{thick, orange}
    \torus{-5}{-4}{2.3}{thick, orange}
    \trefoil{-8}{1.5}{1}{red,very thick}
    \hopflink{10.5}{0}{1}{red,very thick}
    \unknot{0.5}{-5}{.75}{red,very thick}
    \unknot{3}{-4.5}{.75}{red,very thick}
    \trefoil{3}{9}{1.5}{red,very thick}
    \hopflink{-5}{11}{1.5}{red,very thick}
    \unknot{9}{10}{1.5}{red,very thick}
    \node at (-8,-.5) {\tiny$v_j$};
    \node at (1.6,-6) {\tiny$v_k$};
    \node at (10.5,-1.5) {\tiny$v_l$};
    \node at (-12,-2.5) {\tiny${\color{red}L_j}$};
\node at (-2,-8.5) {\tiny${\color{red}L_k}$};
\node at (7.5,-3.5) {\tiny${\color{red}L_l}$};
\node at (-2,13) {\tiny${\color{red}L}$};
\node at (-1.5,6) {\tiny${\color{orange}\Sigma}$};
\node at (1.5,1) {\tiny${\color{orange}\Sigma}$};
    \end{tikzpicture}
}

\newcommand{\lasagnafillingrelfigure}[1]{
\begin{tikzpicture}[anchorbase,scale=#1]
    \begin{scope}[shift={(-18,-1)},xscale=3,yscale=12]
        \draw[dashed] (0,1) to [out=0,in=90] (1,0) to [out=270,in=0] (0,-1)
        (0,1) to [out=180,in=90] (-1,0) to [out=270,in=180] (0,-1);
    \end{scope}
    \begin{scope}[shift={(-18,-1)},xscale=9,yscale=12]
        \draw[thick] (0,-1) to [out=15,in=180] (2,-1.3) to [out=0,in=270] (3.8,0)
        to [out=90,in=0] (2,1.3) to [out=180,in=340] (0,1);
    \end{scope}
    \bgsphere{-4}{-4}{10}
    \fgsphere{-4}{-4}{10}
    \sphere{1}{-5}{3}
    \sphere{10.5}{-.5}{3}
    \sphere{-8}{.6}{3}
    \draw[white, line width=1mm]
    (-9.5,.7) to [out=80,in=290] (-7.3,11)
    (-2.9,11) to [out=280, in=180] (-1,9) to [out=0,in=260] (.7,11)
    (5.2,10.7) to [out=250, in=180] (6,8) to [out=0,in=260] (7.3,10)
    (10.3,10) to [out=250, in=90] (9,6) to [out=270,in=110] (12,0)
    (9,0) to [out=115, in=0] (0,5) to [out=180,in=70] (-6.5,.7)
    (3.1,-4.7) to [out=75, in=0] (-.5,3) to [out=180,in=110] (-.8,-5.2)
    (1.5,-4.4) to [out=75, in=0] (-.1,.5) to [out=180,in=110] (.7,-5.2)
    ;
    \draw[orange, thick]
    (-9.5,.7) to [out=80,in=290] (-7.3,11)
    (-2.9,11) to [out=280, in=180] (-1,9) to [out=0,in=260] (.7,11)
    (5.2,10.7) to [out=250, in=180] (6,8) to [out=0,in=260] (7.3,10)
    (10.3,10) to [out=250, in=90] (9,6) to [out=270,in=110] (12,0)
    (9,0) to [out=115, in=0] (0,5) to [out=180,in=70] (-6.5,.7)
    (3.1,-4.7) to [out=75, in=0] (-.5,3) to [out=180,in=110] (-.8,-5.2)
    (1.5,-4.4) to [out=75, in=0] (-.1,.5) to [out=180,in=110] (.7,-5.2)
    ;
    \torushole{2}{7}{2}{thick, orange}
    \torushole{-.3}{1.5}{2}{thick, orange}
    \torus{-5}{-4}{2.3}{thick, orange}
    \trefoil{-8}{1.5}{1}{red,very thick}
    \hopflink{10.5}{0}{1}{red,very thick}
    \unknot{0.5}{-5}{.75}{red,very thick}
    \unknot{3}{-4.5}{.75}{red,very thick}
    \begin{scope}[shift={(-6.9,4.5)},xscale=2.3,yscale=.4]
        \draw[red, very thick, opacity=.5]
        (-.5,0) to [out=90,in=180] (.5,1) to [out=0,in=90] (1.5,0) to [out=270,in=0] (.5,-1) to [out=180,in=270] (-.5,0);
    \end{scope}
    \trefoil{3}{9}{1.5}{red,very thick}
    \hopflink{-5}{11}{1.5}{red,very thick}
    \unknot{9}{10}{1.5}{red,very thick}
    \node at (-8,-8) {\small$F$};
    \node at (-8,-.5) {\tiny$v_j$};
    \node at (1.6,-6) {\tiny$v_k$};
    \node at (10.5,-1.5) {\tiny$v_l$};
    \end{tikzpicture}
\quad \sim \quad
\begin{tikzpicture}[anchorbase,scale=#1]
    \begin{scope}[shift={(-18,-1)},xscale=3,yscale=12]
        \draw[dashed] (0,1) to [out=0,in=90] (1,0) to [out=270,in=0] (0,-1)
        (0,1) to [out=180,in=90] (-1,0) to [out=270,in=180] (0,-1);
    \end{scope}
    \begin{scope}[shift={(-18,-1)},xscale=9,yscale=12]
        \draw[thick] (0,-1) to [out=15,in=180] (2,-1.3) to [out=0,in=270] (3.8,0)
        to [out=90,in=0] (2,1.3) to [out=180,in=340] (0,1);
    \end{scope}
    \bsphere{-4}{-4}{10}
    \fsphere{-4}{-4}{10}
    \sphere{10.5}{-.5}{3}
    \draw[white, line width=1mm]
    (-8,4.5) to [out=80,in=290] (-7.3,11)
    (-2.9,11) to [out=280, in=180] (-1,9) to [out=0,in=260] (.7,11)
    (5.2,10.7) to [out=250, in=180] (6,8) to [out=0,in=260] (7.3,10)
    (10.3,10) to [out=250, in=90] (9,6) to [out=270,in=110] (12,0)
    (9,0) to [out=115, in=0] (0,5) to [out=180,in=20] (-3.3,4.5)
    ;
    \draw[orange, thick]
    (-8,4.5) to [out=80,in=290] (-7.3,11)
    (-2.9,11) to [out=280, in=180] (-1,9) to [out=0,in=260] (.7,11)
    (5.2,10.7) to [out=250, in=180] (6,8) to [out=0,in=260] (7.3,10)
    (10.3,10) to [out=250, in=90] (9,6) to [out=270,in=110] (12,0)
    (9,0) to [out=115, in=0] (0,5) to [out=180,in=20] (-3.3,4.5)
    ;
    \torushole{2}{7}{2}{thick, orange}
    \hopflink{10.5}{0}{1}{red,very thick}
    \begin{scope}[shift={(-6.9,4.5)},xscale=2.3,yscale=.4]
        \draw[red, very thick]
        (-.5,0) to [out=90,in=180] (.5,1) to [out=0,in=90] (1.5,0) to [out=270,in=0] (.5,-1) to [out=180,in=270] (-.5,0);
    \end{scope}
    \trefoil{3}{9}{1.5}{red,very thick}
    \hopflink{-5}{11}{1.5}{red,very thick}
    \unknot{9}{10}{1.5}{red,very thick}
    \node at (-8,3.4) {\tiny$v_i$};
    \node at (10.5,-1.5) {\tiny$v_l$};
    \end{tikzpicture}
    }

%% file: maintext.tex
\setcounter{tocdepth}{1}
\tableofcontents

\section{Introduction}

The study of link homology theories has revealed profound connections between
low-dimensional topology, representation theory, and higher category theory.
Originally pioneered by Khovanov~\cite{Kho} in the form of a categorification of
the Jones polynomial, link homology provides powerful knot invariants
that detect subtle topological and geometric structures. By a link homology
theory I mean a functor
\begin{equation*}
\linkhom \colon \Links(\R^3) \longrightarrow \Kb(R\gpmod),
\end{equation*}
from the category of links in $\R^3$ and link cobordisms in $\R^3\times I$ to the bounded homotopy
category of chain complexes of graded projective $R$-modules, where $R$ denotes a graded
commutative ring. The most important examples in what follows are the \emph{general linear
link homologies}, pioneered by Khovanov and
Rozansky~\cite{KR}, which categorify the Reshetikhin--Turaev invariants
for $\glN$, see \S\ref{sec:glnlinkhom}.

This survey reviews recent advances connecting link homology theories to
invariants of smooth 4-manifolds\footnote{All \fourm-manifolds considered in this text
are compact, smooth and oriented and links are framed and oriented, unless stated
otherwise.} and extended topological quantum field theories (TQFTs). The
starting point is joint work with Morrison and Walker~\cite{MWW}, which showed
that link homology theories satisfying certain functoriality and monoidality
requirements extend far beyond invariants of links in $\R^3$. They give rise to
skein-theoretic invariants of oriented smooth $4$-manifolds with boundary links,
now known as \emph{skein lasagna modules}\footnote{The name comes from
\cite{MN}. \emph{Lasagna diagrams} appear in \cite{MWW} as higher-dimensional
analog of the \emph{spaghetti and meatball pictures} for planar algebras
attributed to Jones \cite{walker06}.}. These invariants admit computation along
handle decomposition \cite{MN,MWW2,HRW3}, furnish invariants of embedded and
immersed surfaces \cite{MWW3}, and exhibit striking sensitivity---they
distinguish smooth structures on $4$-manifolds~\cite{ren2024khovanov} and detect
exotic surfaces \cite{sullivan2025barnatanskeinlasagnamodules}, while vanishing
on $\mathbb{CP}^2$ and $S^2\times S^2$ \cite{sullivan2024kirby}. 

The framework of skein theory also clarifies how link homology fits into
the broader context of extended topological quantum field theories, as organized
by the \emph{cobordism hypothesis}, the \emph{tangle hypothesis}, and the
\emph{periodic table} of $k$-tuply monoidal
$(n-k)$-categories~\cite{MR1355899,MR2555928}. Classical skein modules of
$3$-manifolds, first studied by Przytycki and Turaev~\cite{MR1194712,MR964255},
form a $(3+\epsilon)$-dimensional TQFT\footnote{An
$(n+\epsilon)$-dimensional TQFT is an extended TQFT defined on manifolds of
dimension between $0$ and $n$, with diffeomorphisms of $n$-manifolds acting as
isomorphisms on the top level.}, which is completely determined by its value on
the point\footnote{The project of making this statement rigorous is still ongoing, see e.g.~\cite{Scheimbauer-thesis,MR4228258,MR4536120}.}, namely the ribbon category encoding the local relations of the
underlying link invariants, e.g. the Jones polynomial. Skein lasagna
modules likewise form a $(4+\epsilon)$-dimensional TQFT,
determined by a braided monoidal $2$-category that encodes the local relations
of the underlying link homology theory, e.g. Khovanov homology. 

Skein theory based on link homology thus realizes part of the vision of
Crane--Frenkel~\cite{MR1295461}: to use \emph{categorification}, as motivated by
Lusztig's canonical bases~\cite{MR1035415}, to construct algebraically
computable 4-dimensional TQFTs of sensitivity comparable to Donaldson invariants
(and later: Seiberg--Witten invariants). A fascinating early perspective on this
proposal and its connection to many contemporary developments, which proved
important, such as Kapranov--Voevodsky's braided monoidal
$2$-categories~\cite{MR1278735}, the movie moves of Carter--Saito~\cite{MR1238875},
and the Crane--Yetter state-sum invariants for 4-manifolds~\cite{CraneYetter},
appears in Baez's \emph{This Week's Finds} \cite{BaezWeeksFinds}. 

\paragraph{Outline} This survey, written on the occasion of the 2025
International Congress of Basic Science, aims to provide an accessible overview
of the state of the art in topological quantum field theories based on link
homology via skein theory. After reviewing the extension from link homology to
skein theory (Theorem~\ref{thm:skeinfromLH}), I discuss the general linear link
homologies $\KhRN$ and the proof of their functoriality in $S^3$ via the
sweep-around move, followed by the operadic setting for the definition of skein
lasagna modules and the relation to braided monoidal $2$-categories. Basic properties,
computational techniques, and applications are collected in \S\ref{sec:properties-computations}, with a focus on invariants of
embedded and immersed surfaces, handle-attachment formulas, and explicit
computations leading to the detection of exotic smooth structures in
\cite{ren2024khovanov}. Finally, \S\ref{sec:TQFTcontext} situates skein
lasagna modules in the context of extended TQFTs and discusses recent progress
toward homotopy-coherent, chain-level versions.


\section{Skein Theory from Link Homology}

The following theorem summarizes the main constructions of \cite{MWW} and is an extension of \cite[Theorem 2.1]{MWW3}, allowing as target category $\target$ an arbitrary symmetric
monoidal cocomplete $1$-category, whose tensor product preserves colimits
separately in each variable, and adding conclusion 4. A typical example is
$\target = \Der(R\gmod)$, the derived category of chain complexes of graded
modules over a (graded) commutative ring $R$, with link homology factoring
through $\Kb(R\gpmod)$, the bounded homotopy category of graded projective $R$-modules.

\begin{thm}\label{thm:skeinfromLH} 
Suppose we are given a \emph{link homology} functor for links in $\R^3$,
\[
\linkhom \colon \Links(\R^3) \longrightarrow \target,
\]
such that $\linkhom$ is
\renewcommand{\labelenumi}{\alph{enumi}.}
\begin{enumerate}
    \item invariant under the trace of $2\pi$ rotation of $\R^3$.
    \item monoidal\footnote{Lax monoidality is sufficient. This structure makes $\linkhom$ compatible with the $\mathbb{E}_3$-monoidal structures on $\Links(\R^3)$ and $\target$, which are symmetric as we are working with plain 1-categories.} under disjoint union: $\linkhom(L_1)\otimes \linkhom(L_2) \to \linkhom(L_1 \sqcup L_2)$,
    \item invariant under the sweep-around move, \cite[(1.1)]{MWW}:
\end{enumerate}
    \begin{equation}
        \label{eq:sweep}
    \includegraphics[width=.7\textwidth]{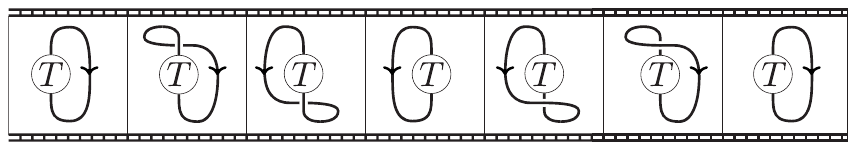}
    \end{equation}

\noindent Then $\linkhom$ extends naturally to:
\renewcommand{\labelenumi}{\arabic{enumi}.}
\begin{enumerate}
    \item A link homology for links in $S^3$ \cite[Definition 4.8]{MWW}, with values in $\target$.
    \item An algebra for the lasagna operad \cite[\S5.1]{MWW}, with values in $\target$.
    \item A $(4+\epsilon)$-dimensional TQFT, whose top-dimensional layer is given by skein modules associated to pairs $(W,L)$ of $4$-manifolds $W$ with links $L \subset \partial W$ taking values in $\target$ \cite[\S5.2]{MWW}.
    \item A locally $\target$-enriched braided monoidal 2-category with duals for objects, adjoints for 1-morphisms, and the analog of a ribbon structure~\cite[\S6]{MWW}.
\end{enumerate}
\end{thm}

The constructions summarized in Theorem~\ref{thm:skeinfromLH}, which I survey in the following subsections, were already implicitly present in unpublished work of Morrison and Walker around 2007 \cite{walker06, walker07}, and in particular influenced the formulation of blob homology \cite{MR2806651,MR2978449}.  
The theorem becomes meaningful once a link homology theory $\linkhom$ satisfying the listed properties is provided. 

\subsection{General linear link homology}
\label{sec:glnlinkhom}
A main contribution of \cite{MWW} is to show that the $\glN$ link homology
theories satisfy the hypotheses of Theorem~\ref{thm:skeinfromLH}. These link
homology theories are categorifications of the Reshetikhin--Turaev link
invariants associated with the complex Lie algebra $\glN$. They were pioneered
by Khovanov and Rozansky \cite{KR} and have been rediscovered and reconstructed
using a variety of mathematical techniques, see e.g,
\cite{stroppel2022categorification} for a recent survey. Especially useful for our purposes is the
combinatorial formulation using foams \cite{MSV,LQR,QR,RW,RoW}, which is
functorial for links in $\R^3$ by \cite{ETW} following \cite{Bla} for the case
$N=2$:

    \begin{equation}
        \label{eq:CKhRR3}
    \CKhRN \colon \Links(\R^3) \longrightarrow \Kb(\Z\gpmod),
    \end{equation}

The definition of the functor $\CKhRN$ proceeds diagrammatically and assigns algebraic data at three levels:
\begin{itemize}
    \item \textbf{Link diagrams:} to each generic planar projection of a link, i.e. to each \emph{link diagram}, $\CKhRN$ assigns a chain complex.
    \item \textbf{Movies of link diagrams:} each link cobordism, visualized as a movie of link diagrams, is decomposed into elementary movies (Reidemeister moves and Morse moves), and $\CKhRN$ assigns a chain map to each movie.
    \item \textbf{Movie moves:} the assignments must respect isotopies
    of cobordism relative to the boundary; specifically, for movies related by the so-called
    Carter--Saito movie moves~\cite{MR1238875}, the associated chain maps must be homotopic.
\end{itemize}

Taking the homology of chain complexes produced by \eqref{eq:CKhRR3} yields a link homology theory $\KhRN$ valued in
the category of bigraded abelian groups
\begin{equation}
\label{eq:KhRR3}
\KhRN \colon \Links(\R^3) \longrightarrow \bigrabgrps
\end{equation}
where one grading is by internal \emph{quantum degree}, and
the other is the \emph{homological grading} of the chain complex. A framed link
cobordism $\Sigma$ from $L_0$ to $L_1$ induces a well-defined homomorphism
\[
\KhRN(\Sigma) \colon \KhRN(L_0) \to \KhRN(L_1),
\]
that is homogeneous of quantum degree $(1-N)\chi(\Sigma)$ and homological degree
zero. 

Moreover, the requirements a. and b. of Theorem~\ref{thm:skeinfromLH} are
straightforward to verify for $\CKhRN$ and $\KhRN$ by means of the diagrammatic
construction.

\begin{exa}
The $\mathfrak{gl}_2$ homology agrees with Khovanov homology \cite{Kho} up to changes in normalization and, possibly, passing to the mirror link.
\end{exa}

\begin{exa}
The $\mathfrak{gl}_1$ homology of any framed link $L$ is free of rank $1$, supported
in quantum degree $-f$ and homological degree $f$, where $f$ denotes the self-linking of $L$ \cite[\S 3.2]{MWW3}.
\end{exa}

\subsection{Link homology in the 3-sphere}
The $3$-sphere $S^3$ may be regarded as a one-point compactification of $\R^3$. A link
    in $S^3$ generically avoids the point at infinity and can thus be assumed to live in $\R^3$. Link cobordisms in $S^3 \times [0,1]$ can similarly be modeled entirely
    within $\R^3 \times [0,1]$, away from infinity. The only subtlety arises when
    considering isotopies of cobordisms that pass through the point at infinity,
    which give rise to the \emph{sweep-around move} \eqref{eq:sweep}.

\begin{thm}
    \label{thm:sweep} Let $N\in \Z_{\geq 1}$, then the $\glN$ link homology functor \eqref{eq:CKhRR3}
    assigns the identity morphism to every instance of the sweep-around move
    \eqref{eq:sweep} and hence satisfies the hypotheses of
    Theorem~\ref{thm:skeinfromLH}. 
\end{thm}

The difficulty of proving the sweep-around move, the main reason for the delay
between \cite{walker07} and \cite{MWW}, is that it is not a single move but an
infinite family of moves, indexed by choices of tangles $T$ in \eqref{eq:sweep}, which
are non-local, at least from the perspective of link diagrams, far away from the
point at infinity. The key idea for our proof in \cite{MWW} is a
categorification of the \emph{Kauffman trick} from \cite[Lemma 2.4]{MR899057},
which exploits the dependence between the skein relation and Reidemeister moves
of type 2 and 3. In our case, this enables a systematic comparison of chain maps
associated to Reidemeister moves of type 3, which happen when the closure strand
in the sweep-around move passes either in front or past the back of a diagram of
the tangle $T$. This technique of proof applies to variations of $\glN$ homology
\cite{MWW3,ren2024khovanov,sullivan2025barnatanskeinlasagnamodules} and has been
adapted to related settings \cite{MR4504654,chen2025flipmapinvolutionskhovanov}.

\smallskip
Once the sweep-around move is established, the extension of the link homology functor to links in $S^3$, as asserted in
Theorem~\ref{thm:skeinfromLH}.1, proceeds in two main steps. Here I discuss
them only at the level of links and refer to \cite[\S 4.1-2]{MWW} for details on
the behavior under link cobordisms.

\begin{itemize}
\item \textbf{Removal of the parametrization of $\R^3$.} The groupoid of parametrizations of $\R^3$ forms a torsor over the group of orientation-preserving diffeomorphisms $\mathrm{Diff}^+(\R^3)$, which is path-connected with fundamental group 
$\pi_1(\mathrm{Diff}^+(\R^3)) \cong \Z/2\Z$, generated by a $2\pi$ rotation of $\R^3$. 
By functoriality of the link homology in $\R^3$, any smooth path of parametrizations induces a link isotopy, and hence an isomorphism on link homology. 
Assumption a. in Theorem~\ref{thm:skeinfromLH} guarantees that these isomorphisms depend only on the endpoints of the path, not on the particular choice of isotopy. 
The resulting invariant of a link in the unparametrized ambient $\R^3$ is therefore the \emph{transitive system} of all such homologies equipped with the canonical isomorphisms; equivalently, it can be described as the (co)limit over the groupoid of parametrizations.

\item \textbf{Extension to links in unparametrized 3-spheres.} To define the
invariant for a link $L$ in $S^3$, observe that any choice of base point $p \in
S^3 \setminus L$ presents $L$ as a link in the $3$-ball $S^3 \setminus \{p\}$.
Moving $p$ along a path in the link complement induces a canonical isomorphism
between the corresponding link homologies. Since the fundamental group of the
complement is generated by meridians of $L$, it suffices to check monodromy
around these loops, which is precisely captured by the \emph{sweep-around move}.
Assumption c. in Theorem~\ref{thm:skeinfromLH} ensures that this monodromy acts
trivially. The invariant of $L$ in $S^3$ is defined as the
\emph{transitive system} of all such link homologies equipped with the canonical
isomorphisms, equivalently described as the (co)limit over the fundamental
groupoid of the complement $S^3 \setminus L$. 
\end{itemize}

Recent work of Teng \cite{teng2025endkhovanovhomologyexotic} uses functoriality
under link cobordisms in $S^3\times [0,1]$ to define a version of Khovanov
homology capable of detecting exotic Lagrangian and symplectic planes in $\R^4$.

\subsection{Lasagna algebra}

We are now ready to consider the relevant types of skeins. From now on we let $W$ be a \fourm-manifold and $L
\subset \partial W$ a link, unless stated otherwise. We also fix $N\in \Z_{\geq 1}$ and work with the link homology theory $\KhRN$ from \S\ref{sec:glnlinkhom}.

\begin{definition} One defines:

\begin{enumerate} 
    \item A \emph{lasagna skein} $F = (\Sigma, \{(B_i,L_i)\})$ of $W$ with boundary $L$ consists of:
\begin{itemize}
    \item A finite collection of disjointly embedded $4$-balls $B_i\hookrightarrow \operatorname{int}W$; and
    \item A framed oriented surface $\Sigma$ properly embedded in $W \setminus \cup_i \operatorname{int} B_i$, meeting $\partial W$ in $L$ and each $\partial B_i$ in a link $L_i$.
\end{itemize}
\item A \emph{lasagna filling} of $W$ with boundary $L$ is a lasagna skein as above with:
\begin{itemize}
        \item For each $i$, a homogeneous label $v_i \in \KhRN(\partial B_i, L_i)$. 
\end{itemize}
\item A \emph{lasagna diagram} is a lasagna skein for $W=B^4$.
\item The \emph{lasagna operad} is a colored operad with \begin{itemize}
\item set of colors given by the framed oriented links in $S^3$;
\item set of operations given by lasagna diagrams as above, with the $L_i$ serving as inputs and $L$ as output; and
\item composition given by gluing a lasagna diagram to an input sphere of another lasagna diagram.
\end{itemize}
\end{enumerate}
\end{definition}

\begin{figure}[ht]
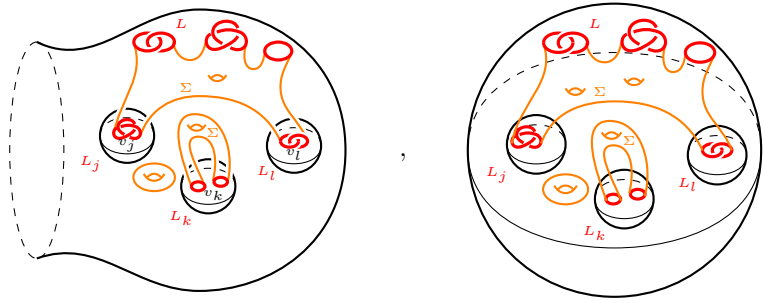

    \[
    \lasagnafillingfigure{.12}\qquad,\qquad \lasagnafigure{.13}
    \]
    \caption{A lasagna filling of a generic \fourm-manifold $W$ and a lasagna diagram.}
\end{figure}
\noindent For a comparison with
concepts from skein theory for 3-manifolds, see Table~\ref{tab:skein
lasagna-comparison}.

\begin{theorem}[{\cite[Theorem 5.2]{MWW}}]\label{thm:lasagnaalgebra} The link homology functor $\KhRN$ from \eqref{eq:KhRR3} extends to an algebra for the lasagna operad.
\end{theorem}
 This means that $\KhRN$ assigns a bigraded abelian group $\KhRN(L)$ to each link $L\subset S^3$ and further, a homogeneous morphism 
\[
\KhRN(\Sigma) \colon \bigotimes_i \KhRN(\partial B_i, L_i) \to \KhRN(\partial W, L),
\]
to every lasagna diagram, such that gluing lasagna diagrams is compatible with
composing morphisms. The idea of the proof uses that the relevant maps are
already provided by the functoriality statement of Theorem~\ref{thm:skeinfromLH}.1 in the
case of at most one input link. When considering more than one input link, the input
spheres first have to be tubed together along an embedded graph, yielding a
single input sphere with a split disjoint union of input links. One then uses
the (lax) monoidality to define the associated map. Independence on the choice
of tubing graph is a consequence of the sweep-around move. 

\begin{exa}
    Each lasagna filling $F$ of $B^4$ with boundary $L$ and surface $\Sigma$ yields an element $\KhRN(F) \in \KhRN(L)$ by evaluating $\KhRN(\Sigma)$ on the tensor product of input labels $v_i$ of $F$.
\end{exa}

\subsection{Skein modules for 4-manifolds}

An algebra for the lasagna operad provides both the labeling data for skein
modules of \fourm-manifolds as well as the skein relations. As before we let $W$
be a \fourm-manifold, $L \subset \partial W$ a link, and $N\in \Z_{\geq 1}$. I will describe three
equivalent definitions of the skein modules associated to the link homology $\KhRN$.

\begin{definition}
    The \emph{skein lasagna module} 
\[
\Sz(W; L) := \Z\langle \text{lasagna fillings $F$ of $W$ with boundary $L$}\rangle/\sim,
\]
is obtained as quotient of the bigraded abelian group freely generated by all lasagna fillings of $W$ with boundary $L$, by the subgroup that enforce the
transitive and linear closure of the following relations:
\begin{enumerate}
    \item Linearity in the labels $v_i\in \KhRN(L_i)$.
    \item Equivalence under replacement of an input ball $B_i$ with a lasagna filling $F$ of a $4$-ball such that $v_i = \KhRN(F)$, followed by isotopy rel $\partial W$:
\[
\lasagnafillingrelfigure{.13}
\]
\end{enumerate}
\end{definition}

For the following alternative description, let $\mathcal{C}(W,L)$ denote the set of lasagna skeins of $W$ with boundary $L$. For each such skein $S\in \mathcal{C}(W,L)$, we let $\mathrm{input}(S)$ denote the finite set of input links $L_i$ in $S$.

\begin{definition}
    The \emph{skein lasagna module} is the quotient
\[
\Sz(W; L) = \bigoplus_{S\in \mathcal{C}(W,L)} \bigotimes_{L_i\in \mathrm{input}(S)} \KhRN(L_i)/\sim,
\]
by the subgroup generated by the relators
\[
\underbrace{(\otimes_{k\in K} v_k) \otimes \KhRN(\Sigma)(\otimes_{j\in J} v_j)}_{\text{in } \bigotimes_{L_i\in \mathrm{input}(S)} \KhRN(L_i)}
-
\underbrace{(\otimes_{k\in K} v_k) \otimes (\otimes_{j\in J} v_j)}_{\text{in } \bigotimes_{L_i\in \mathrm{input}(S')}\KhRN(L_i)}
\]  Here $S,S'$ are lasagna
skeins such that $S'$ is obtained from $S$ by attaching a lasagna diagram, i.e.
a lasagna skein of $B^4$, with underlying surface $\Sigma$, followed by an
isotopy. The set $J$ indexes input links of $\Sigma$, there is a unique output
link of $\Sigma$, which also serves as input link for $S$ and $K$ indexes the
remaining input links of $S$, which then also appear in $S'$.
\end{definition}

The following reformulation appears in \cite{ren2024khovanov}. To formulate it,
one considers $\mathcal{C}(W,L)$ as a category with morphisms $S' \to S$ generated by attachments of lasagna diagrams to input spheres of $S$, yielding $S'$ up to a specified isotopy rel boundary, and with composition given by iterated attachments and composed isotopies. As a consequence of Theorem~\ref{thm:lasagnaalgebra}, one obtains a functor:
\begin{align*}
    \KhRN\colon \mathcal{C}(W,L) &\to \bigrabgrps, \qquad S \mapsto \textstyle \bigotimes_{L_i\in \mathrm{input}(S)} \KhRN(L_i)
\end{align*}

Here we suppress from the notation certain grading shifts, which depend on the lasagna skeins and have the effect of making all morphisms in the image grading-preserving.

\begin{definition}
    The \emph{skein lasagna module} is the colimit 
\[
\Sz(W; L) = \operatorname{colim}\left( \mathcal{C}(W,L) \xrightarrow{\KhRN} \bigrabgrps\right),
\]
\end{definition}

Note that this description directly generalizes to link homology theories with
values in other symmetric monoidal cocomplete target categories $\target$,
provided the conditions of Theorem~\ref{thm:skeinfromLH} are satisfied.

\begin{rem} The setting of Theorem~\ref{thm:skeinfromLH} is not the only one in which an extension from a categorical link invariant to a skein-theoretic \fourm-manifold invariant can be envisioned. For example, a construction based on link Floer homology appears in \cite{chen2022floerlasagnamoduleslink}. It uses a modified notion of skeins that accommodates links with multiple basepoints and link cobordisms with embedded arcs that connect basepoints.
\end{rem}

\subsection{Braided monoidal 2-categories}
\label{sec:braided-2-cat}

Given a link homology theory $\linkhom$ satisfying the hypotheses of
Theorem~\ref{thm:skeinfromLH}, the construction summarized there produces, in
addition to a functorial link homology in $S^3$, a lasagna algebra, and skein
lasagna modules, a \emph{locally $\target$-enriched braided monoidal 2-category}
$\mathcal{C}_{\linkhom}$. This construction was first provided in \cite[\S
6]{MWW} for the prototypical case of $\KhRN$, carefully accounting for
(semi-)strictness; here I describe the generalization only informally.

Objects of $\mathcal{C}_{\linkhom}$ correspond to finite sequences of framed
points in a 2-disk $D^2\subset \R^2$, 1-morphisms are tangles in $D^2 \times
[0,1]$ with horizontal composition implemented by stacking. The 2-morphisms
between two tangles $S$ and $T$ with the same boundary are computed as the link
homology $H(T\cup_{\partial S= \partial T} \overline{S})\in \target$ of the link
obtained by gluing $T$ to the mirror-reverse of $S$. The horizontal and vertical
composition of 2-morphisms is implemented using the functoriality of $\linkhom$ under
link cobordisms, leading to $\target$-enriched hom categories. The braided
monoidal structure on this 2-category is inherited from the naturality of the
construction under embedding little 2-disks in larger 2-disks. The requisite
braiding data \cite{MR1278735, MR1402727} on the level of 1-morphisms and 2-morphisms can be described
explicitly in terms of certain shuffle braids and tangle cobordisms
respectively. By construction, all objects in $\mathcal{C}_{\linkhom}$ admit
duals and all 1-morphisms admit adjoints \cite[\S 6.4]{MWW} and a categorification
of the ribbon equation holds \cite[\S 6.5]{MWW}.

Conceptually, this braided monoidal 2-category plays a role for the link homology theory $\linkhom$ and its skein lasagna modules for smooth 4-manifolds analogous to that played by ribbon categories, e.g. of quantum group representations, in the decategorified setting of Reshetikhin--Turaev invariants \cite{MR1036112} and associated skein modules of 3-manifolds, see Table~\ref{tab:skein lasagna-comparison}. We will return to its role in determining extended TQFTs in \S\ref{sec:TQFTcontext}.

\begin{table}[h]
\centering
\begin{tabular}{|p{3.8cm}|p{3.4cm}|p{4.0cm}|}
\hline
\textbf{feature} & \textbf{skein modules} & \textbf{skein lasagna modules} \\
 \hline
link invariant & Reshetikhin--Turaev  & link homology theory $\linkhom$ \\
 \hline
categorical data & ribbon category  & ribbon 2-category $\mathcal{C}_\linkhom$ \\
\hline
ambient manifold & oriented 3d  & oriented smooth 4d \\
\hline
type of skeins & ribbon graphs & lasagna fillings \\
\hline
local labelling at& coupons & input balls \\
\hline
boundary condition&  points &  links \\
\hline
labelling by & morphisms & link homology classes\\
\hline
\end{tabular}
\caption{Comparison of features between classical skein modules for 3-manifolds based on ribbon categories and skein lasagna modules for surfaces in 4-manifolds based on a link homology theory.}
\label{tab:skein lasagna-comparison}
\end{table}

\section{Properties, Computations, and Applications}
\label{sec:properties-computations}

This section discusses structural properties, computational techniques, and
selected applications of $\glN$ skein lasagna modules. 
\subsection{Basic Properties}

This subsection collects the foundational properties of $\glN$ skein lasagna modules that arise straightforwardly from their construction. These include functoriality, monoidality under disjoint union, and the behavior under standard gluing operations.

\begin{itemize}
    \item \textbf{Recovery of link homology:} For the 4-ball $B^4$ and a link $L \subset \partial B^4 \cong S^3$, we have a canonical isomorphism
    \[
    \KhRN(L) \xrightarrow{\cong} \Sz(B^4; L).
    \]
    induced by decorating the radial skein. This further illustrates why the extension of $\KhRN$ to links in $S^3$ is essential for the skein module construction.

    \item \textbf{Gradings:} The skein module is $\Z\times \Z$-graded by quantum and homological degree and decomposes further according to classes in relative second homology:
    \[
    \Sz(W; L) = \bigoplus_{\alpha \in H_2^L(W;\Z)} \Sz(W; L, \alpha).
    \]
    Here $H_2^L(W; \Z) := \partial^{-1}([L]) \subset H_2(W, L; \Z)$ is the preimage of the fundamental class of $L$ under the connecting map $\partial$ of the long exact sequence for relative homology; it is a torsor over $H_2(W; \Z)$.  

    \item \textbf{Gluing and functoriality under inclusions:} When a 4-manifold $W = W_1 \cup_Y W_2$ is obtained by gluing 4-manifolds $W_1, W_2$ along a common part $Y$ of their boundaries, this induces a map on skein modules
    \begin{equation}
        \label{eq:gluing}
    \Sz(W_1; L_1\cup L_Y) \otimes \Sz(W_2; \overline{L_Y}\cup L_2) \xrightarrow{~} \Sz(W; L_1\cup L_2).
    \end{equation}
    Here the boundary links are presented as unions of tangles $L_1,L_2, L_Y,
    \overline{L_Y}$ and $L_Y$ gets glued to $\overline{L_Y}$.  As a consequence,
    skein modules are functorial under embeddings of \fourm-manifolds that are
    compatible with the boundary links. A detailed discussion appears in
    \cite[\S2.2]{MWW}. More generally, \eqref{eq:gluing} can be upgraded to a
    presentation of $\Sz(W; L_1\cup L_2)$ as relative tensor product of modules
    $\Sz(W_1; L_1\cup -)$ and $\Sz(W_2; \overline{-}\cup L_2)$ over a suitably
    defined linear skein category associated to $(Y; \partial{L_1})$ \cite[4.4.2]{walker06}.
    \item \textbf{Monoidality under disjoint unions and (boundary) connected sum:} The skein module is (laxly) monoidal under disjoint union and over a field $\kk$, this is a strong monoidal equivalence:
    \[
   \Sz(W_1; L_1;\kk) \otimes \Sz(W_2; L_2;\kk)\xrightarrow{\cong}  \Sz(W_1 \sqcup W_2; L_1 \sqcup L_2;\kk),
    \]
    The skein module also behaves monoidally under both connected sum and boundary connected sum \cite[\S 7]{MN}.
\end{itemize}

\subsection{Handle Attachments}
\label{subsec:handle-attachments}
A powerful strategy for computing skein lasagna modules, introduced by Manolescu and Neithalath~\cite{MN} and fully developed in~\cite{MWW2}, is to proceed inductively along a handle decomposition of the $4$-manifold $W$. Given a link $L \subset \partial W$ and a handle decomposition of $W$ ordered by index, the skein module of $(W; L)$ is computed \emph{in reverse}, by successively removing handles and analyzing the their effect on the skein module—possibly altering the boundary link in the process. The computation terminates in a disjoint union of $4$-balls, where the skein module reduces to a link homology calculation.

I now briefly summarize the effect of removing handles of each index on the $\glN$ skein lasagna module, following the reverse computation strategy. For simplicity we work over a field.\vspace{-1mm}
\begin{itemize}
\item \textbf{Four-handles:} Attaching a 4-handle induces an isomorphism on skein modules, so 4-handles can also be freely removed, c.f.~\cite[Proposition 2.1]{MN}. 
    
    \item \textbf{Three-handles:} Attaching a 3-handle induces a surjection on skein modules, also proven in~\cite[Proposition 2.1]{MN}. The kernel of this surjection is described in~\cite[\S 3.2]{MWW2} as image of the difference of cobordism maps associated with the two attaching hemispheres of the 3-handle.

\item \textbf{Two-handles:} The \emph{Manolescu--Neithalath 2-handle formula} shows that the effect of attaching a 2-handle can simulated by inserting parallel cables of the attaching knot, performing a symmetrization procedure, and assembling the results into a filtered colimit as the number of strands tends to infinity. This approach was developed in~\cite{MN, MWW2} and further perspectives in terms of Kirby-colored link homology are discussed in~\cite{HRW3,vonmerkl2025computingcoloredkhovanovhomology}.

\item \textbf{One-handles:} Every 1-handle corresponds to a boundary connected sum---possibly a self-sum---along disjoint 3-balls in the boundary, and typically interacts with the boundary link. Algebraically, the effect on skein modules is described by computing co-invariants for a skein category associated to $B^3$, which acts on both components of the attaching region of the 1-handle~\cite[\S 4]{MWW2}.

\item \textbf{Zero-handles:} After all higher-index handles have been removed, the manifold becomes a disjoint union of $4$-balls. The skein module is then computed as a tensor product of link homologies over the remaining boundary links.
\end{itemize}

In both the 1-handle and 2-handle cases, the reduction to the skein module of a
manifold with the handle detached comes at the cost of considering infinite
families boundary links. For 2-handles, these are parallel cables of the
attaching knot, forming a natural and often partially computable family. In
contrast, the 1-handle formula is significantly more difficult to control: the
set of resulting boundary links is much less structured, and the co-invariant
constructions involve actions by categories that are not yet well understood.
As shown in~\cite[Theorem 1.5]{MWW2}, this can lead to skein modules that are
not locally finite-dimensional.

It is an open question, for which 4-manifolds $W$ and links $L$ the skein module
$\Sz(W; L)$ is of finite-rank in each $\Z\times \Z \times H_2^L(W;\Z)$-degree.
Work in preparation by Qi--Robert--Sussan--Wagner addresses a refined question
of finite generation by considering additional symmetries on equivariant skein
lasagna modules.

\subsection{Invariants of embedded surfaces}
Skein lasagna modules provide a natural home for invariants of smoothly
embedded---and even immersed---surfaces in 4-manifolds, just as skein modules based
on ribbon categories serve as targets for invariants of framed links in
3-manifolds.

A key subtlety is that skein lasagna modules are spanned by lasagna fillings,
which consist of oriented, framed surfaces. Since not every embedded surface
admits a framing, one instead works with punctured surfaces: finitely many 4-balls
are excised from the ambient manifold, puncturing the surface so that what remains becomes framable. The new boundary components can then be canonically
decorated with specific link homology classes corresponding to nonzero-framed
unknots, arising from cobordism maps associated to the Reidemeister move of type 1.
This construction extends further to singular surfaces. Isolated
singularities---such as transverse double points---can be modeled by removing
neighborhoods and decorating their links (e.g., Hopf links) with canonical
homology classes. In this way, skein lasagna modules yield invariants of
immersed as well as embedded surfaces.

These skein elements can be viewed as generalizations of the \emph{relative
Khovanov--Jacobsson classes} associated to surfaces in the 4-ball
\cite{Jac,MR4562563}, which are known to distinguish certain exotic pairs of
surfaces---embedded surfaces that are topologically but not smoothly isotopic
\cite{MR4726569}.

Such skein classes also form the basis for extracting topological information about
smooth surfaces in 4-manifolds, for instance lower bounds on the minimal genus
in a given relative second homology class. In \cite{MWW3}, $\GLN$-equivariant
$\glN$ link homology is used to establish such bounds by analyzing the grading
support of skein modules modulo torsion, over the base ring
$H^*(\mathrm{B}\GLN)$. A related approach by Ren and Willis
\cite{ren2024khovanov} uses a Lee-type deformation of $\mathfrak{gl}_2$ link
homology to construct a filtered skein module, with the quantum filtration
yielding lower bounds on the genus function. In the case of the 4-ball, this
recovers Rasmussen's s-invariant \cite{Ras, MR2462446}. A different approach to
extending the s-invariants to surfaces in certain other 4-manifold appears in
\cite{MR4541332}. 

Analyzing the torsion in Bar-Natan type deformations of $\mathfrak{gl}_2$ skein
modules enables the detection of exotic pairs of knotted surfaces that remain
exotic after one internal stabilization
\cite{sullivan2025barnatanskeinlasagnamodules}.

\subsection{Computations and Sensitivity towards Exotica}
A range of explicit computations of skein lasagna modules have been carried out
for small 4-manifolds and certain classes of links in the boundary, with a view towards testing the sensitivity of the invariant. 

\begin{table}[h]
\centering
\begin{tabular}{|c|c|c|c|}
\hline
$W^4$ & $L$ & $\Sztwo(W^4;L)$ & \textbf{Reference} \\
\hline
$S^4$ & $\emptyset$ & $\kk$ & \cite{MWW} \\
\hline
$B^3 \times S^1$  & $\sqcup_{2m} S^1$ & loc. finite rank exactly for $m\leq 1$ & \cite{MWW2} \\
\hline
$S^2 \times D^2$ & $\emptyset$ & $\kk[x,x^{-1}]/x\kk[x]$ per $H_2$ class & \cite{MN} \\
 & L & $\Sztwo(S^2 \times D^2;\emptyset) \otimes \mathrm{RW}(L)$ & \cite{sullivan2024kirby} \\
\hline
$S^2 \times S^2$ & $\emptyset$ & 0 & \cite{sullivan2024kirby} \\
\hline
$\mathbb{CP}^2$ & $\emptyset$ & $0$ & \cite{MN,ren2024khovanov} \\
\hline
$\overline{\mathbb{CP}}^2$ & $\emptyset$ & nonzero, description conjectural& \cite{MN,ren2024khovanov,vonmerkl2025computingcoloredkhovanovhomology} \\
\hline
\end{tabular}
\caption{Sample computations of skein lasagna modules for various 4-manifolds.}
\label{tab:skein-computations}
\end{table}
The sensitivity with respect to orientation was demonstrated in \cite{MN} through
partial computations for $\mathbb{CP}^2$ with both orientations. For links $L$
in $S^2\times S^1$, the invariant depends on the bulk 4-manifold: $B^3\times S^1$ often leads to invariants that are non locally
finite-rank \cite{MWW2}, while $S^2\times D^2$ produces locally finite rank
tensor multiples of the Rozansky--Willis invariant $\mathrm{RW}(L)$
\cite{rozansky2010categorification,MR4332675}. The latter description, together
with the vanishing for $S^2 \times S^2$, was established by
Sullivan--Zhang~\cite{sullivan2024kirby}.

The remarkable recent preprint of Ren and Willis \cite{ren2024khovanov} contains
many additional computations beyond Table~\ref{tab:skein-computations},
including a comparison of the skein modules of the exotic pair of knot traces
$X_{{-}1}(-5_2)$ and $X_{{-}1}(P(3,-3,8))$. In generating second homology
classes and homological degree zero, these skein modules differ in quantum
degree $-1$ over $\Q$. Thus skein lasagna modules can detect exotic smooth
structure by purely algebro-combinatorial means. Also notable are vanishing
results for skein lasagna modules, e.g. for $4$-manifolds that contain a positive
self-intersection embedded $S^2$ \cite[Theorem 1.4]{ren2024khovanov}, which use the fact that having vanishing skein module is a property
that is inherited under embeddings.

\section{Topological Quantum Field Theory Context}
\label{sec:TQFTcontext}

The skein lasagna modules described in this survey arise as the 4-dimensional
layer of an extended local topological quantum field theory (TQFT), whose associated
manifold invariants can be described in skein-theoretic terms for oriented
manifolds of dimensions up to $4$. This theory is determined by the braided
monoidal 2-category $\mathcal{C}_{\linkhom}$ extracted in
\S\ref{sec:braided-2-cat} from the underlying link homology $\linkhom$. 

Braided monoidal 2-categories appear in the periodic table of $n$-categories
\cite{MR1355899} as categorification of braided monoidal
categories\footnote{The colors in the table indicate equal total categorical dimension
$n=(n-k)+k$; categorification usually means passing from a cell to the adjacent cell with higher column index.}. Link homology and functoriality under
cobordism maps are natural consequences, in analogy to how Reshetkhin--Turaev
tangle invariants are captured by ribbon categories.

\begin{center}
	\begin{tabular}{|c||c| c| c| c| c|} 
		\hline
		 $\mathbb{E}_k \backslash${\small $n-k$}& 0 & 1 & 2 & $\cdots$ \\ [0.5ex] 
		\hline\hline
		$-$ & \cellcolor{blue!25} sets & \cellcolor{green!25} categories &  \cellcolor{yellow!25} 2-categories & $\cdots$\\ 
		\hline
		$\mathbb{E}_1$ & \cellcolor{green!25} monoids & \cellcolor{yellow!25} monoidal cats & \cellcolor{orange!25} monoidal 2-cats&$\cdots$ \\
		\hline
		$\mathbb{E}_2$ &\cellcolor{gray!25} comm. monoids & \cellcolor{orange!25} braided cats& \cellcolor{red!25} braided 2-cats&$\cdots$\\
		\hline
		$\mathbb{E}_3$ &\cellcolor{gray!25} $\qquad$ --- $\qquad$ &\cellcolor{gray!25} sym. mon. cats & \cellcolor{purple!25} sylleptic 2-cats &$\cdots$\\
		\hline
		$\mathbb{E}_4$ &\cellcolor{gray!25} --- &\cellcolor{gray!25} --- &\cellcolor{gray!25} sym. mon. 2-cats &$\cdots$\\
		\hline
		$\vdots$ &\cellcolor{gray!25} $\qquad\;\;\;$---$\qquad\;\;\;$ &\cellcolor{gray!25} $\qquad\;\;\;$---$\qquad\;\;\;$ &\cellcolor{gray!25} $\qquad\quad\;$---$\qquad\quad\;$& $\ddots$\\ [1ex] 
		\hline
	   \end{tabular}
  \end{center}

This perspective situates skein lasagna modules within a broader hierarchy of
skein-theoretic topological quantum field theories. Classical instances arise
from the 2-dimensional graphical calculus of monoidal categories or the
3-dimensional graphical calculus of braided monoidal categories. These are often
referred to as the \emph{Turaev--Viro} and \emph{Crane--Yetter} families of TQFTs,
though strictly speaking those names are more closely tied to the corresponding
state-sum models, which extend such theories up by one
dimension, provided additional strong finiteness conditions are satisfied~\cite{TuraevViro, CraneYetter}.

In my Frontiers of Science Award Lecture at the International Congress of Basic Sciences 2025~\cite{WedICBStalk}, I outlined the current progress toward categorified analogues of these theories based on monoidal 2-categories and braided monoidal 2-categories, emphasizing in particular the role of the chosen target higher category.

\begin{center}
	\begin{tabular}{|c||c| c| c| c| c|} 
		\hline		 
		&  linear & loc. linear  & loc. stable  \\ 
		 $\mathbb{E}_k \backslash${\small $n-k$}&  categories & 2-categories & $(\infty,2)$-categories \\ [0.5ex] 
		\hline\hline
		monoidal & \cellcolor{yellow!25} Turaev--Viro 
    & \cellcolor{orange!25} Asaeda--Frohman--Kaiser 
    &  \cite{HRW4} \ref{sec:HRW} \cellcolor{orange!25} \\
		\hline
		braided & \cellcolor{orange!25} Crane--Yetter 
    & \cellcolor{red!25}  \cite{MWW} 
    &  \cite{liu2024braided} \ref{sec:LMGRSW} \cellcolor{red!25}\\ 
		\hline
	   \end{tabular}
  \end{center}

  The main difference is whether one considers local enrichment in a symmetric
monoidal 1-category, such as (graded) abelian groups or vector spaces, or in
fact in a symmetric monoidal stable $\infty$-category, such as chain complexes.
In the former case, monoidal 2-categories lead to skein theories
studied by Asaeda--Frohman and Kaiser~\cite{MR2370224,kaiser2022barnatan} and
many others, or to the invariants of Douglas--Reutter~\cite{douglas2018fusion},
while braided monoidal 2-categories lead to lasagna skein modules as in
\cite{MWW}. 

The next two sections outline the state of chain level versions, which arise when the base of enrichment is upgraded to a symmetric monoidal stable $\infty$-category, such as chain complexes, i.e.~the last column of the above table.

\subsection{Towards a chain level version}
\label{sec:LMGRSW}

Link homology theories are typically constructed from chain complexes associated
to link diagrams. This naturally raises the question whether one can promote the
functoriality of link homology to the chain level,
in such a way that it becomes homotopy-coherent. 

\begin{conjecture}
    \label{conj:chainlinkh} The chain level $\glN$ link invariant from \eqref{eq:CKhRR3} arises as truncation of an $\mathbb{E}_3$-monoidal functor of $(\infty,1)$-categories 
    \begin{equation}
        \label{eq:CKhRR4}
    \CKhRN \colon \Links_\infty(\R^3) \longrightarrow \mathrm{Ch}^b(\Z\gpmod),
    \end{equation}
from the $(\infty,1)$-category $\Links_\infty(\R^3)$ of links in $\R^3$, link
cobordisms in $\R^3 \times I$, isotopies, and higher isotopies, to the
$(\infty,1)$-category of bounded chain complexes, chain maps, homotopies, and higher homotopies.
\end{conjecture}

The traditional approach to functoriality via movie moves breaks down at the
chain level, since one encounters an infinite hierarchy of higher relations
between movies. Instead, a more conceptual framework is required. The guiding
idea is to capture the necessary homotopy-coherence through local data: namely,
an $\mathbb{E}_2$-monoidal $(\infty,2)$-category, generated by $2$-dualizable
objects and equipped with an $\mathrm{SO}(4)$-homotopy-fixed-point structure.
This would serve as the chain-level, categorified analogue of the ribbon
categories underlying the Reshetikhin--Turaev tangle invariants, with the
functor \eqref{eq:CKhRR4} recoverable by restriction to tangles without
endpoints.

While a complete chain-level theory is still under development, several
precursor results point in this direction. In joint work with Stroppel (2021),
we established the homotopy-coherent naturality of the Rouquier
braiding~\cite{Rou2Braid} in a concrete dg model for chain complexes of Soergel
bimodules. This result, now documented
in~\cite{stroppel2024braidingtypesoergelbimodules}, led to Conjecture 3.8 in
Stroppel's ICM article~\cite{stroppel2022categorification}.

In collaboration with Liu, Mazel-Gee, Reutter, and
Stroppel~\cite{liu2024braided}, we resolved this conjecture by constructing an
$\mathbb{E}_2$-monoidal $(\infty,2)$-category of chain complexes of Soergel
bimodules. This structure underlies braid invariants feeding into $\glN$ link
homologies and provides a categorical foundation for triply graded homology
theories~\cite{MR2339573}. Further joint work with Dyckerhoff~\cite{DW25}
relates the braiding on complexes of Soergel bimodules to the concept of perverse
schobers, as proposed by Kapranov--Schechtman~\cite{KS14a,KS14b,KS21}. A key
ingredient in this connection are singular Soergel
bimodules~\cite{Williamson-thesis}, which can be obtained from Soergel bimodules via
higher-categorical idempotent completion~\cite{R25}.

A current limitation of the $\mathbb{E}_2$-monoidal $(\infty,2)$-category of
Soergel complexes is that its generating object is not dualizable. As a
consequence, the associated invariant extends to braids but not to tangles or
general tangle cobordisms, leaving Conjecture~\ref{conj:chainlinkh} unresolved. We plan to remedy this issue by proceeding to $\glN$ quotients in future work.

Homotopy-coherent, chain-level link invariants would have important advantages
beyond their intrinsic structural appeal. They have the potential to capture
finer topological information and improve the computability of associated
invariants. This is already visible for Rouquier complexes of braids: when considering braid closures in the annulus via categorical
traces~\cite{2002.06110}, higher homotopical data is essential for cabling
operations~\cite[\S 6.4]{2019arXiv190404481G}. More generally, chain-level
tangle invariants provide the natural framework for categorifying skein
algebras, avoiding the semisimplification procedures otherwise
required \cite{1806.03416}. 

\subsection{Towards a chain level version in 3d}
\label{sec:HRW}

The guiding idea for a chain-level version of skein theory for surfaces in
$3$-manifolds is that it should yield categorified, partially defined analogs of
TQFTs in the Turaev--Viro family. Conceptually, such a theory ought to be based
on an $\mathbb{E}_1$-monoidal locally stable $(\infty,2)$-category generated by
2-dualizable objects, equipped with an $\mathrm{SO}(3)$-homotopy-fixed-point
structure. Assuming link homology has been modeled locally and homotopy
coherently as outlined in \S\ref{sec:LMGRSW}, such a structure can be
obtained by forgetting from $\mathbb{E}_2$-monoidality to $\mathbb{E}_1$, i.e.\
by discarding the braiding. Since the braiding is the most intricate part of the
higher-categorical data, $\mathbb{E}_1$-examples are comparatively more
accessible directly---for instance by viewing the locally linear monoidal $2$-categories
underlying the Asaeda--Frohman--Kaiser TQFT as enriched in chain complexes.

This approach was initiated in~\cite{HRW4}, which begins with the Bar--Natan
monoidal $2$-category, the categorification of the Temperley--Lieb monoidal category
underlying combinatorial constructions of Khovanov homology \cite{BN2}. Mirroring
Roberts's skein-theoretic description of the Turaev--Viro
theory~\cite{MR1362787}, we obtain a categorified analog of the
Turaev--Viro~TQFT in low dimensions. The main result of~\cite{HRW4} is the
explicit construction and characterization of the invariant of $2$-dimensional
$1$-handlebodies, which takes the form of certain dg categories in the simplest
cases. These categories are generated by objects parametrized by spin networks
adapted to a triangulation of the surface. Their hom pairings yield power series
in the variable $q$, whose graded Euler characteristics recover the
Turaev--Viro hermitian pairing when $q$ is specialized to a complex root of
unity. Ongoing work with Hogancamp and Rose extends this construction to the
$3$-dimensional level and closed surfaces. 

\paragraph{Acknowledgements}
I am deeply grateful to Kim Morrison and Kevin Walker for our collaboration and
many illuminating discussions, and to the organizers of the 2025 ICBS for the
opportunity to present this work.\footnote{Portions of this article have been edited for clarity using generative AI.}

\paragraph{Funding}
I acknowledge support from the Deutsche
Forschungsgemeinschaft (DFG, German Research Foundation) under Germany's
Excellence Strategy - EXC 2121 ``Quantum Universe'' - 390833306 and the
Collaborative Research Center - SFB 1624 ``Higher structures, moduli spaces and
integrability'' - 506632645.

%% file: ICBS_Wedrich-arxiv.bbl
\begin{thebibliography}{LMGRSW24}
    
\bibitem[AF07]{MR2370224}
Marta Asaeda and Charles Frohman.
\newblock A note on the {B}ar-{N}atan skein module.
\newblock {\em Internat. J. Math.}, 18(10):1225--1243, 2007.

\bibitem[Bae]{BaezWeeksFinds}
John~C. Baez.
\newblock This {W}eek's {F}inds.
\newblock available at
  {\href{https://math.ucr.edu/home/baez/twf.html}{https://math.ucr.edu/home/baez/twf.html}},
  accessed 25.08.2025.

\bibitem[BD95]{MR1355899}
John~C. Baez and James Dolan.
\newblock Higher-dimensional algebra and topological quantum field theory.
\newblock {\em J. Math. Phys.}, 36(11):6073--6105, 1995.

\bibitem[BJS21]{MR4228258}
Adrien Brochier, David Jordan, and Noah Snyder.
\newblock On dualizability of braided tensor categories.
\newblock {\em Compos. Math.}, 157(3):435--483, 2021.

\bibitem[Bla10]{Bla}
Christian Blanchet.
\newblock An oriented model for {K}hovanov homology.
\newblock {\em J. Knot Theory Ramifications}, 19(2):291--312, 2010.

\bibitem[BN96]{MR1402727}
John~C. Baez and Martin Neuchl.
\newblock Higher-dimensional algebra. {I}. {B}raided monoidal {$2$}-categories.
\newblock {\em Adv. Math.}, 121(2):196--244, 1996.

\bibitem[BN05]{BN2}
Dror Bar-Natan.
\newblock Khovanov's homology for tangles and cobordisms.
\newblock {\em Geom. Topol.}, 9:1443--1499, 2005.

\bibitem[BW08]{MR2462446}
Anna Beliakova and Stephan Wehrli.
\newblock Categorification of the colored {J}ones polynomial and {R}asmussen
  invariant of links.
\newblock {\em Canad. J. Math.}, 60(6):1240--1266, 2008.

\bibitem[CF94]{MR1295461}
Louis Crane and Igor~B. Frenkel.
\newblock Four-dimensional topological quantum field theory, {H}opf categories,
  and the canonical bases.
\newblock volume~35, pages 5136--5154. 1994.

\bibitem[Che22]{chen2022floerlasagnamoduleslink}
Daren Chen.
\newblock Floer lasagna modules from link {F}loer homology, 2022.
\newblock \arxiv{2203.07650}.

\bibitem[Coo23]{MR4536120}
Juliet Cooke.
\newblock Excision of skein categories and factorisation homology.
\newblock {\em Adv. Math.}, 414:Paper No. 108848, 51, 2023.

\bibitem[CS93]{MR1238875}
J.~Scott Carter and Masahico Saito.
\newblock Reidemeister moves for surface isotopies and their interpretation as
  moves to movies.
\newblock {\em J. Knot Theory Ramifications}, 2(3):251--284, 1993.

\bibitem[CY93]{CraneYetter}
Louis Crane and David Yetter.
\newblock A categorical construction of {$4$}d topological quantum field
  theories.
\newblock In {\em Quantum topology}, volume~3 of {\em Ser. Knots Everything},
  pages 120--130. World Sci. Publ., River Edge, NJ, 1993.

\bibitem[CY25]{chen2025flipmapinvolutionskhovanov}
Daren Chen and Hongjian Yang.
\newblock The flip map and involutions on Khovanov homology, 2025.
\newblock \arxiv{2506.00824}.

\bibitem[DR18]{douglas2018fusion}
Christopher~L. Douglas and David~J. Reutter.
\newblock Fusion 2-categories and a state-sum invariant for 4-manifolds, 2018.
\newblock \arxiv{1812.11933}.

\bibitem[DW25]{DW25}
Tobias Dyckerhoff and Paul Wedrich.
\newblock Perverse schobers of coxeter type $\mathbb{A}$, 2025.
\newblock \arxiv{2504.08496}.

\bibitem[ETW18]{ETW}
Michael Ehrig, Daniel Tubbenhauer, and Paul Wedrich.
\newblock Functoriality of colored link homologies.
\newblock {\em Proc. Lond. Math. Soc. (3)}, 117(5):996--1040, 2018.

\bibitem[GHW21]{2002.06110}
Eugene Gorsky, Matthew Hogancamp, and Paul Wedrich.
\newblock Derived traces of {S}oergel categories.
\newblock {\em Int. Math. Res. Not. IMRN}, 2021.

\bibitem[GW23]{2019arXiv190404481G}
Eugene Gorsky and Paul Wedrich.
\newblock Evaluations of annular {K}hovanov-{R}ozansky homology.
\newblock {\em Math. Z.}, 303(1):Paper No. 25, 57, 2023.

\bibitem[HRW22]{HRW3}
Matthew {Hogancamp}, D.~E.~V. {Rose}, and Paul {Wedrich}.
\newblock A {K}irby color for {K}hovanov homology, 2022.
\newblock \arxiv{2210.05640}, to appear in \textit{J. Eur. Math. Soc.}

\bibitem[HRW24]{HRW4}
Matthew {Hogancamp}, D.~E.~V. {Rose}, and Paul {Wedrich}.
\newblock Bordered invariants from {K}hovanov homology, 2024.
\newblock \arxiv{2404.06301}.

\bibitem[HS24]{MR4726569}
Kyle Hayden and Isaac Sundberg.
\newblock Khovanov homology and exotic surfaces in the 4-ball.
\newblock {\em J. Reine Angew. Math.}, 809:217--246, 2024.

\bibitem[Jac04]{Jac}
M.~Jacobsson.
\newblock An invariant of link cobordisms from {K}hovanov homology.
\newblock {\em Algebr. Geom. Topol.}, 4:1211--1251 (electronic), 2004.

\bibitem[Kai25]{kaiser2022barnatan}
Uwe Kaiser.
\newblock Bar-{N}atan theory and tunneling between incompressible surfaces in
  3-manifolds.
\newblock {\em Topology Appl.}, 369:Paper No. 109390, 54, 2025.

\bibitem[Kau87]{MR899057}
Louis~H. Kauffman.
\newblock State models and the {J}ones polynomial.
\newblock {\em Topology}, 26(3):395--407, 1987.

\bibitem[Kho00]{Kho}
M.~Khovanov.
\newblock A categorification of the {J}ones polynomial.
\newblock {\em Duke Math. J.}, 101(3):359--426, 2000.

\bibitem[Kho07]{MR2339573}
Mikhail Khovanov.
\newblock Triply-graded link homology and {H}ochschild homology of {S}oergel
  bimodules.
\newblock {\em Internat. J. Math.}, 18(8):869--885, 2007.

\bibitem[KR08]{KR}
Mikhail Khovanov and Lev Rozansky.
\newblock Matrix factorizations and link homology.
\newblock {\em Fund. Math.}, 199(1):1--91, 2008.

\bibitem[KS15]{KS14b}
Mikhail Kapranov and Vadim Schechtman.
\newblock Perverse schobers, 2015.
\newblock \arxiv{1411.2772}.

\bibitem[KS16]{KS14a}
Mikhail Kapranov and Vadim Schechtman.
\newblock Perverse sheaves over real hyperplane arrangements.
\newblock {\em Ann. of Math. (2)}, 183(2):619--679, 2016.

\bibitem[KS25]{KS21}
Mikhail Kapranov and Vadim Schechtman.
\newblock P{ROB}s and perverse sheaves {I}: symmetric products.
\newblock {\em Selecta Math. (N.S.)}, 31(2):Paper No. 23, 32, 2025.

\bibitem[KV94]{MR1278735}
M.~M. Kapranov and V.~A. Voevodsky.
\newblock {$2$}-categories and {Z}amolodchikov tetrahedra equations.
\newblock In {\em Algebraic groups and their generalizations: quantum and
  infinite-dimensional methods ({U}niversity {P}ark, {PA}, 1991)}, volume~56 of
  {\em Proc. Sympos. Pure Math.}, pages 177--259. Amer. Math. Soc., Providence,
  RI, 1994.

\bibitem[LMGRSW24]{liu2024braided}
Yu~Leon Liu, Aaron Mazel-Gee, David Reutter, Catharina Stroppel, and Paul
  Wedrich.
\newblock A braided $(\infty,2)$-category of {S}oergel bimodules, 2024.
\newblock \arxiv{2401.02956}.

\bibitem[LQR15]{LQR}
A.D. Lauda, H.~Queffelec, and D.E.V. Rose.
\newblock Khovanov homology is a skew {H}owe $2$-representation of categorified
  quantum $\mathfrak{sl}(m)$.
\newblock {\em Algebr. Geom. Topol.}, 15(5):2517--2608, 2015.

\bibitem[LS22]{MR4504654}
Robert Lipshitz and Sucharit Sarkar.
\newblock A mixed invariant of nonorientable surfaces in equivariant {K}hovanov
  homology.
\newblock {\em Trans. Amer. Math. Soc.}, 375(12):8807--8849, 2022.

\bibitem[Lur09]{MR2555928}
Jacob Lurie.
\newblock On the classification of topological field theories.
\newblock In {\em Current developments in mathematics, 2008}, pages 129--280.
  Int. Press, Somerville, MA, 2009.

\bibitem[Lus90]{MR1035415}
G.~Lusztig.
\newblock Canonical bases arising from quantized enveloping algebras.
\newblock {\em J. Amer. Math. Soc.}, 3(2):447--498, 1990.

\bibitem[MMSW23]{MR4541332}
Ciprian Manolescu, Marco Marengon, Sucharit Sarkar, and Michael Willis.
\newblock A generalization of {R}asmussen's invariant, with applications to
  surfaces in some four-manifolds.
\newblock {\em Duke Math. J.}, 172(2):231--311, 2023.

\bibitem[MN22]{MN}
Ciprian Manolescu and Ikshu Neithalath.
\newblock Skein lasagna modules for 2-handlebodies.
\newblock {\em J. Reine Angew. Math.}, 788:37--76, 2022.

\bibitem[MSV09]{MSV}
M.~Mackaay, M.~Sto{\v{s}}i{\'c}, and P.~Vaz.
\newblock {$\mathfrak{sl}_{N}$}-link homology {$(N\geq 4)$} using foams and the
  {K}apustin--{L}i formula.
\newblock {\em Geom. Topol.}, 13(2):1075--1128, 2009.

\bibitem[MW11]{MR2806651}
Scott Morrison and Kevin Walker.
\newblock Higher categories, colimits, and the blob complex.
\newblock {\em Proc. Natl. Acad. Sci. USA}, 108(20):8139--8145, 2011.

\bibitem[MW12]{MR2978449}
Scott Morrison and Kevin Walker.
\newblock Blob homology.
\newblock {\em Geom. Topol.}, 16(3):1481--1607, 2012.

\bibitem[MWW22]{MWW}
Scott Morrison, Kevin Walker, and Paul Wedrich.
\newblock Invariants of 4-manifolds from {K}hovanov-{R}ozansky link homology.
\newblock {\em Geom. Topol.}, 26(8):3367--3420, 2022.

\bibitem[MWW23]{MWW2}
Ciprian Manolescu, Kevin Walker, and Paul Wedrich.
\newblock Skein lasagna modules and handle decompositions.
\newblock {\em Adv. Math.}, 425:Paper No. 109071, 40, 2023.

\bibitem[MWW24]{MWW3}
Scott Morrison, Kevin Walker, and Paul Wedrich.
\newblock Invariants of surfaces in smooth 4-manifolds from link homology,
  2024.
\newblock \arxiv{2401.06600}.

\bibitem[Prz91]{MR1194712}
J\'{o}zef~H. Przytycki.
\newblock Skein modules of {$3$}-manifolds.
\newblock {\em Bull. Polish Acad. Sci. Math.}, 39(1-2):91--100, 1991.

\bibitem[QR16]{QR}
H.~Queffelec and D.E.V. Rose.
\newblock The $\mathfrak{sl}_n$ foam $2$-category: a combinatorial formulation
  of {K}hovanov--{R}ozansky homology via categorical skew {H}owe duality.
\newblock {\em Adv. Math.}, 302:1251--1339, 2016.

\bibitem[QW21]{1806.03416}
Hoel {Queffelec} and Paul {Wedrich}.
\newblock {Khovanov homology and categorification of skein modules}.
\newblock {\em Quantum Topol.}, 12(1):129--209, 2021.

\bibitem[Ras10]{Ras}
J.~Rasmussen.
\newblock Khovanov homology and the slice genus.
\newblock {\em Invent. Math.}, 182(2):419--447, 2010.

\bibitem[Rec25]{R25}
Isabela Recio.
\newblock Higher idempotent completion for Soergel bimodules, 2025.
\newblock \arxiv{2508.00767}.

\bibitem[Rob95]{MR1362787}
Justin Roberts.
\newblock Skein theory and {T}uraev-{V}iro invariants.
\newblock {\em Topology}, 34(4):771--787, 1995.

\bibitem[Rou17]{Rou2Braid}
R.~Rouquier.
\newblock {K}hovanov-{R}ozansky homology and $2$-braid groups.
\newblock In {\em Categorification in geometry, topology, and physics}, 2017.

\bibitem[Roz10]{rozansky2010categorification}
Lev Rozansky.
\newblock A categorification of the stable {SU}(2)
  {W}itten-{R}eshetikhin-{T}uraev invariant of links in {S2 x S1}, 2010.
\newblock \arxiv{1011.1958}.

\bibitem[RT90]{MR1036112}
N.~Yu. Reshetikhin and V.~G. Turaev.
\newblock Ribbon graphs and their invariants derived from quantum groups.
\newblock {\em Comm. Math. Phys.}, 127(1):1--26, 1990.

\bibitem[RW16]{RW}
David~E.V. Rose and Paul Wedrich.
\newblock Deformations of colored $\mathfrak{sl}(n)$ link homologies via foams.
\newblock {\em Geom. Topol.}, 20(6):3431--3517, 2016.

\bibitem[RW20]{RoW}
Louis-Hadrien Robert and Emmanuel Wagner.
\newblock A closed formula for the evaluation of foams.
\newblock {\em Quantum Topol.}, 11(3):411--487, 2020.

\bibitem[RW24]{ren2024khovanov}
Qiuyu Ren and Michael Willis.
\newblock Khovanov homology and exotic $4$-manifolds, 2024.
\newblock \arxiv{2402.10452}.

\bibitem[{Sch}14]{Scheimbauer-thesis}
Claudia {Scheimbauer}.
\newblock {Factorization homology as a fully extended topological field
  theory}, 2014.
\newblock Dissertation, ETH Zurich,
  \url{http://scheimbauer.at/ScheimbauerThesis.pdf}.

\bibitem[SS22]{MR4562563}
Isaac Sundberg and Jonah Swann.
\newblock Relative {K}hovanov-{J}acobsson classes.
\newblock {\em Algebr. Geom. Topol.}, 22(8):3983--4008, 2022.

\bibitem[Str23]{stroppel2022categorification}
Catharina Stroppel.
\newblock Categorification: tangle invariants and {TQFT}s.
\newblock In {\em I{CM}---{I}nternational {C}ongress of {M}athematicians.
  {V}ol. {II}. {P}lenary lectures}, pages 1312--1353. EMS Press, Berlin, [2023]
  \copyright 2023.

\bibitem[Sul25]{sullivan2025barnatanskeinlasagnamodules}
Ian~A. Sullivan.
\newblock Bar-{N}atan skein lasagna modules and exotic surfaces in 4-manifolds,
  2025.
\newblock \arxiv{2504.03968}.

\bibitem[SW24]{stroppel2024braidingtypesoergelbimodules}
Catharina Stroppel and Paul Wedrich.
\newblock Braiding on type {A} {S}oergel bimodules: semistrictness and
  naturality, 2024.
\newblock \arxiv{2412.20587}.

\bibitem[SZ24]{sullivan2024kirby}
Ian~A. Sullivan and Melissa Zhang.
\newblock Kirby belts, categorified projectors, and the skein lasagna module of
  ${S}^{2}\times{{S}^{2}}$, 2024.
\newblock \arxiv{2402.01081}.

\bibitem[Ten25]{teng2025endkhovanovhomologyexotic}
Yikai Teng.
\newblock End Khovanov homology and exotic Lagrangian planes, 2025.
\newblock \arxiv{2510.01151}.

\bibitem[Tur88]{MR964255}
V.~G. Turaev.
\newblock The {C}onway and {K}auffman modules of a solid torus.
\newblock {\em Zap. Nauchn. Sem. Leningrad. Otdel. Mat. Inst. Steklov. (LOMI)},
  167(Issled. Topol. 6):79--89, 190, 1988.

\bibitem[TV92]{TuraevViro}
V.~G. Turaev and O.~Ya. Viro.
\newblock State sum invariants of {$3$}-manifolds and quantum {$6j$}-symbols.
\newblock {\em Topology}, 31(4):865--902, 1992.

\bibitem[vM25]{vonmerkl2025computingcoloredkhovanovhomology}
Karim~Ritter von Merkl.
\newblock Computing colored {K}hovanov homology, 2025.
\newblock \arxiv{2505.03916}.

\bibitem[Wal06]{walker06}
Kevin Walker.
\newblock Tqfts, May 2006.
\newblock Notes available at
  {\href{https://canyon23.net/math/}{https://canyon23.net/math/}}, accessed
  25.08.2025.

\bibitem[Wal07]{walker07}
Kevin Walker.
\newblock Khovanov homology as a {TQFT}, April 2007.
\newblock Talk notes available at
  {\href{https://canyon23.net/math/talks/}{https://canyon23.net/math/talks/}},
  accessed 14.08.2025.

\bibitem[Wed25]{WedICBStalk}
Paul Wedrich.
\newblock From link homology to topological quantum field theories, July 2025.
\newblock 2025 Frontiers of Science Award lecture, recording available at
  \url{https://www.youtube.com/watch?v=G4ZHSR1S_oY}, accessed 19.08.2025.

\bibitem[{Wil}08]{Williamson-thesis}
Geordie {Williamson}.
\newblock {Singular {S}oergel bimodules}, 2008.
\newblock Dissertation, Freiburg,
  \url{http://people.mpim-bonn.mpg.de/geordie/GW-thesis.pdf}.

\bibitem[Wil21]{MR4332675}
Michael Willis.
\newblock Khovanov homology for links in {$\#^r(S^2\times S^1)$}.
\newblock {\em Michigan Math. J.}, 70(4):675--748, 2021.

\end{thebibliography}
